\documentclass{amsart}[1996/10/24]
\usepackage{verbatim}
\usepackage{amssymb,amsmath}
\copyrightinfo{2007}{American Mathematical Society}

\newtheorem{lemma}{Lemma} [section]
\newtheorem{thm}[lemma]{Theorem}
\newtheorem{cor}[lemma]{Corollary}

\theoremstyle{remark}
\newtheorem{remark}[lemma]{Remark}

\newcommand{\F}{{\mathbb F}}

\newcommand{\fq}{\F_q}

\DeclareMathOperator{\Sym}{Sym}
\DeclareMathOperator{\Gal}{Gal}
\newcommand{\Z}{{\mathbb Z}}

\begin{document}



\title{Permutation binomials over finite fields}

\author{Ariane M. Masuda}
\address{
School of Mathematics and Statistics,
Carleton University,
1125 Colonel By Drive,
Ottawa, Ontario, Canada K1S 5B6
}

\curraddr{
Department of Mathematics and Statistics,
University of Ottawa, 
585 King Edward Avenue, 
Ottawa, Ontario, Canada K1N 6N5
}
\email{amasuda@uottawa.ca}

\author{\hbox{Michael E. Zieve}}
\address{
Center for Communications Research,
805 Bunn Drive,
Princeton, New Jersey 08540
}
\email{zieve@math.rutgers.edu}
\urladdr{www.math.rutgers.edu/$\sim$zieve/}

\thanks{The authors thank Jeff VanderKam and Daqing Wan for
valuable conversations, and Igor Shparlinski for suggesting the
use of the Brun--Titchmarsh theorem in section 4.}


\date{December 23, 2008}


\begin{abstract}
We prove that if $x^m + ax^n$ permutes the prime field
$\F_p$, where $m>n>0$ and $a\in\F_p^*$, then
$\gcd(m-n,p-1) > \sqrt{p}-1$.  Conversely, we prove that
if $q\ge 4$ and $m>n>0$ are fixed and satisfy
$\gcd(m-n,q-1) > 2q(\log \log q)/\log q$, then there
exist permutation binomials over $\F_q$ of the form
$x^m + ax^n$ if and only if $\gcd(m,n,q-1) = 1$.
\end{abstract}


\maketitle


\section{Introduction}

A polynomial over a finite field is called a \emph{permutation polynomial}
if it permutes the elements of the field.  These polynomials first arose
in work of Betti~\cite{B} and Hermite~\cite{H} as a way to represent
permutations.  A general theory was developed by Hermite~\cite{H} and
Dickson~\cite{D}, with many subsequent developments by Carlitz and others.
The simplest class of nonconstant polynomials are
the monomials $x^m$ with $m>0$, and one easily checks that $x^m$ permutes
$\F_q$ if and only if $m$ is coprime to $q-1$.  However, already for
binomials the situation becomes much more mysterious.  Some examples occurred
in Hermite's work~\cite{H}, and Mathieu~\cite{M}
showed that $x^{p^i} - a x$ permutes $\F_q$ whenever $a$ is not a
$(p^i-1)$-th power in $\F_q$; here $p$ denotes the characteristic of $\F_q$.

A general nonexistence result was proved by Niederreiter and Robinson~\cite{NR}
and improved by Turnwald~\cite{T}:
\begin{thm}
\label{NR}
If $f(x):=x^m+ax^n$ permutes\/ $\F_q$, where $m>n>0$ and $a\in\F_q^*$, then
either $q\le (m-2)^4 + 4m - 4$ or $m=np^i$.
\end{thm}

This result implies that, when $q>m^4$, the only permutation binomials over $\F_q$
are the compositions of Mathieu's examples with permutation monomials.  The
key ingredient in the proof of Theorem~\ref{NR} is Weil's lower bound~\cite{Weil}
for the number of $\F_q$-rational points on the curve $(f(x)-f(y))/(x-y)$.

We do not know whether Theorem~\ref{NR} can be improved in general.  However,
for prime fields it was improved by Wan~\cite{W} and Turnwald~\cite{T}; by using
ingredients from both of their proofs, one can show the following result, which improves
both of their results:
\begin{thm}
\label{WT} If $f(x):=x^m+ax^n$ permutes the prime field\/ $\F_p$,
where $m>n>0$ and $a\in\F_p^*$, then $p-1 \le
(m-1)\cdot\max(n,\gcd(m-n,p-1))$.
\end{thm}

The proofs of Wan and Turnwald rely on a trick due to
Hermite~\cite{H}, which can be viewed as a character sum argument:
they find an integer $\ell$ with $0<\ell<p-1$ such that $f(x)^{\ell}$ mod
$(x^p-x)$ has degree $p-1$.  This implies that
$\sum_{\alpha\in\F_p}f(\alpha)^{\ell}\ne 0$, so $f$ does not permute $\F_p$.
We will prove the following stronger result by exhibiting two integers
$\ell$, of which at least one must have the above property:
\begin{thm}
\label{intro1}
If $f(x):=x^m+ax^n$ permutes the prime field\/ $\F_p$, where $m>n>0$ and $a\in\F_p^*$,
then $\gcd(m-n,p-1)\ge \sqrt{p-(3/4)}-(1/2) > \sqrt{p} - 1$.
\end{thm}
Writing $g:=\gcd(m-n,p-1)$, the conclusion of this result can be restated as
$p-1 \le (g+1)\cdot g$, whereas the conclusion of Theorem~\ref{WT} says that
$p-1\le (m-1)\cdot\max(n,g)$.  Thus, Theorem~\ref{intro1} implies Theorem~\ref{WT}
whenever $g+1\le m-1$, which always holds except in the special case that $n=1$
and $(m-1)\mid (p-1)$.  We emphasize that Theorem~\ref{intro1} is qualitatively
different from all previous results, since it gives a bound on $p$
which depends only on $\gcd(m-n,p-1)$, and not on the degree of $f$.

Both Theorem~\ref{WT} and Theorem~\ref{intro1} yield improvements to
Weil's lower bound for the number of rational points on the curve
$(f(x)-f(y))/(x-y)$ appearing in the proof of Theorem~\ref{NR}. On a
related note, for any polynomial $f$ over $\F_p$ of degree in a
certain range, Voloch~\cite{V} has improved Weil's upper bound for
this same curve.  In a different direction, for hyperelliptic curves
over $\F_p$  one can improve both the upper and lower Weil bound
when the genus is on the order of $\sqrt{p}$, by using Stepanov's
method~\cite{Korobov, Mitkin, Stark, Mitkin2, SV, Baghdadi,
Zannier}. All of these improvements are specific to prime fields. It
would be interesting to understand what are the types of curves for
which one has such improvements to Weil's bounds.

Theorem~\ref{intro1} is not true for nonprime fields; one counterexample
is $x^{10}+3x$ over $\F_{343}$, and we have found several infinite families
of counterexamples, which we will describe in a forthcoming paper.

Returning to prime fields, we suspect that Theorem~\ref{intro1} can be improved.
We checked via computer that, for $p<10^5$, the hypotheses of
Theorem~\ref{intro1} imply that $\gcd(m-n,p-1)>p/(2\log p)$.
It seems likely that this improved result remains true for larger $p$, but we do
not know a proof.  The best we can do is give a heuristic to the effect
that `at random' there would not be any permutation binomials $x^m+ax^n$ over $\F_q$ with
$\gcd(m-n,q-1)<q/(2\log q)$.  Of course, our examples over nonprime fields show
that this heuristic is not always correct, but those examples exhibit
nonrandom features dependent on the subfield structure of $\F_q$, which is in line
with our `at random' notion.

Conversely, following earlier investigations of Hermite~\cite{H} and Brioschi~\cite{Br,Br2},
Carlitz~\cite{C} studied permutation binomials of the form $x^n(x^{(q-1)/2}+a)$.
He showed that there are permutation binomials of this shape (with $n=1$ and $a\in
\F_q^*$)
whenever $q\ge 7$.  He proved a similar result for the form $x(x^{(q-1)/3}+a)$,
and more generally in a paper with Wells~\cite{CW} he proved
\begin{thm}
\label{cwintro} If $d>0$ and $q\equiv 1\pmod{d}$, where $q$ is
sufficiently large compared to $d$, then for each $n>0$ with
$\gcd(n,q-1)=1$ there exists $a\in\F_q^*$ such that $x^n(x^{(q-1)/d}+a)$
permutes\/ $\F_q$.
\end{thm}
The proof of this result is quite remarkable, as it uses the Weil lower bound
on an auxiliary curve to prove the existence of permutation binomials.  This (and a
generalization in~\cite{WL}) is the only known instance of the Weil
bound being used to prove existence of permutation polynomials.  We give a new
proof of a refined version of Theorem~\ref{cwintro}, which allows us to estimate
the number of such $a$'s:
\begin{thm}
\label{intro2} Pick integers $0<n<m$ such that $\gcd(m,n,q-1)=1$,
and suppose $q\ge 4$. If $\gcd(m-n,q-1)>2q(\log\log q)/\log q$, then
there exists $a\in\F_q^*$ such that $x^m+ax^n$ permutes\/ $\F_q$.
Further, letting $T$ denote the number of values $a\in\F_q$ for
which $x^m+ax^n$ permutes\/ $\F_q$, and putting
$r:=(q-1)/\gcd(m-n,q-1)$, we have
\[
\frac{q-2\sqrt{q}+1}{r^{r-1}} - (r-3)\sqrt{q} - 2 \le \frac{T}{(r-1)!} \le
\frac{q+2\sqrt{q}+1}{r^{r-1}} + (r-3)\sqrt{q}.
\]
\end{thm}
We note that the condition $\gcd(m,n,q-1)=1$ is clearly necessary if
$x^m+ax^n$ is to permute $\F_q$. In some special cases, a weaker
estimate for $T$ was derived in a recent paper by
Laigle-Chapuy~\cite{LC}, via methods quite different from ours.

We checked that, for each $q<10^6$, and for every $m>n>0$ satisfying
$\gcd(m,n,q-1)=1$ and $\gcd(m-n,q-1)>2q/\log q$, there exists
$a\in\F_q^*$ such that $x^m+ax^n$ permutes $\F_q$. Combined with our
previously mentioned computer data, this paints a rather clear
picture of permutation binomials over prime fields.

As a final remark, we note that several papers prove results about
the special binomials $x^m+ax$.  In general, if a binomial has a
term of degree coprime to $q-1$, then one can convert it to this
special form by composing with suitable permutation monomials and
reducing mod $(x^q-x)$.  However, there are binomials for which this is
impossible. For instance, $f(x):=x^{26}+17x^3$ permutes $\F_{139}$,
but the degrees of both terms of $f$ have a common factor with
$138$.

Throughout this paper, $\F_q$ is the field of order $q$, and $p$ is the
characteristic of $\F_q$.  In particular, $p$ is always prime.
We prove Theorem~\ref{intro1} in the next section.  Then in
Section~\ref{sec exist} we prove Theorem~\ref{intro2}, and in the
final section we give the heuristic argument mentioned above.
In an appendix we include a proof of Theorem~\ref{WT}.


\section{Nonexistence results}

In this section we prove Theorem~\ref{intro1} in the following form:

\begin{thm}
\label{gen} Suppose $x^n(x^k+a)$ permutes\/ $\F_p$, where $n,k>0$
and $a\in\F_p^*$.  Then $\gcd(k,p-1) \ge \sqrt{p-(3/4)}-(1/2) > \sqrt{p}-1$.
\end{thm}

Our proof relies on Hermite's criterion~\cite{H,D}:

\begin{lemma}
A polynomial $f\in\fq[x]$ is a permutation polynomial if and only if
\begin{enumerate}
\item for each $i$ with $0<i<q-1$, the reduction of
$f^i$ modulo $x^q-x$ has degree less than $q-1$; and
\item $f$ has precisely one root in $\fq$.
\end{enumerate}
\end{lemma}

\begin{proof}[Proof of Theorem~\ref{gen}]
Pick $j>0$ such that $jk\equiv \gcd(k,p-1)\bmod{(p-1)}$ and $\gcd(j,p-1)=1$;
then $x^n(x^k+a)$ permutes $\F_p$ if and only if
$x^{nj}(x^{\gcd(k,p-1)}+a)$ permutes $\F_p$, so we may assume that $k$ divides
$p-1$.  Suppose $f:=x^n(x^k+a)$ permutes $\F_p$, where $k\mid (p-1)$ and
$k<\sqrt{p-(3/4)}-(1/2)$ (and $n,k>0$ and $a\in\F_p^*$).
Then $k^2+k+1<p$.  Let $r$ be the least integer
such that $r\ge (p-1-k)/k^2$.  Then $r<(p-1-k)/k^2+1$, so
\[ kr < (p-1)/k -1 + k = (k-1)(1 - (p-1)/k) + p - 1 \leq p-1.
\]
Also the inequality $k^2+k+1<p$ implies $(p-1-k)/k^2>1$, so $r>1$.

We will apply Hermite's criterion with exponent $kr$.  To this end,
we compute
\[
f^{kr} = x^{nkr}  (x^k + a)^{kr}= x^{nkr} \sum_{i=0}^{kr}
\binom{kr}{i} a^{kr-i} x^{ki}.
\]
\noindent Write $f^{kr} = \sum_{i=0}^{kr} b_i x^{nkr+ki}$,
where $b_i=\binom{kr}{i} a^{kr-i}$.  Since $0<kr<p$ and $p$ is prime,
each $b_i$ is nonzero.  Thus, the degrees of the terms of $f^{kr}$
are
\[
 nkr, nkr+k, nkr+2k,\ldots, nkr+k^2r.
\]
Since $k^2r\geq p-1-k$,  the degrees
include members of every residue class modulo $p-1$ containing
multiples of $k$. In particular, there is a term of degree divisible
by $p-1$; but, since $0<kr<p-1$, Hermite's criterion implies that
$f^{kr}$ cannot have a unique term of degree divisible by $p-1$, so
there must be more than one such term.  Thus,  $nkr \equiv
-E \bmod{(p-1)}$  for some $E$ with $0 \leq E \leq k^2r-(p-1)$.

Likewise, the degrees of the terms of $f^{k(r-1)}$ are
\[
 nk(r-1), nk(r-1)+k, nk(r-1)+2k, \ldots, nk(r-1)+k^2(r-1).
\]
Since  $k^2(r-1) < p-1-k$, these degrees are all in distinct
classes modulo $p-1$, so by Hermite's criterion none of the degrees
can be divisible by $p-1$. Thus, $nk(r-1) \equiv F \bmod{(p-1)}$
for some $F$ with $k \leq F \leq p-1-k-k^2(r-1)$.

Now we have
\[
E(r-1) \equiv -nkr(r-1) \equiv -Fr \bmod{(p-1)},
\]
so $E(r-1)+Fr$ is a multiple of $p-1$.  But
\begin{align*}
0 < kr&\leq E(r-1) + Fr \\
&\leq k^2r(r-1) - (p-1)(r-1) + (p-1)r - kr -
k^2(r-1)r \\
&= p-1-kr < p-1,
\end{align*}
so $E(r-1)+Fr$ lies between consecutive multiples of $p-1$, a contradiction.
\end{proof}

\begin{remark}
The above proof shows that, if
$\gcd(k,p-1)<\sqrt{p-(3/4)}-(1/2)$, then there exists $i$ with $0<i<p-1$
for which the polynomial
$(x^n(x^k+a))^i$ has a unique term of degree divisible by $p-1$,
contradicting our hypothesis that $x^n(x^k+a)$ permutes $\F_p$.
As discussed in the introduction, we suspect that Theorem~\ref{gen}
can be improved substantially.  However, improving our bound by more
than a constant factor will require a new method: if
$\gcd(k,p-1) \ge \sqrt{2p-(7/4)}-(1/2)$, then there is no $i>0$ for which
$(x^n(x^k+a))^i$ has a unique term of degree divisible by $p-1$.
\end{remark}

We now list some consequences of Theorem~\ref{gen}.

\begin{cor}
\label{2p4p} If $x^n(x^k+a)$ permutes\/ $\F_p$, where $n,k>0$ and $a\in\F_p^*$,
then $\gcd(k,p-1)>4$.
\end{cor}

\begin{proof}
When $p>19$, this is an immediate consequence of Theorem~\ref{gen}.
Otherwise, the result can be verified via computer.
\end{proof}

In case either $(p-1)/2$ or $(p-1)/4$ is prime, Corollary~\ref{2p4p}
was conjectured in~\cite{MPW}.
We proved this conjecture in our previous paper~\cite{MZ}, where moreover
we proved that the hypotheses of Corollary~\ref{2p4p}
imply $\gcd(k,p-1)\notin\{2,4\}$ (without assuming primality of $(p-1)/2$
or $(p-1)/4$).  Our proof in~\cite{MZ} did not rely on any computer
calculations; instead we used repeated applications of Hermite's
criterion in several different cases (depending on the class of $p$ mod~$16$).
By using a computer to verify small cases, we can go much further than
Corollary~\ref{2p4p}.  For instance:

\begin{cor}
Suppose $x^n(x^k+a)$ permutes\/ $\F_p$, where $n,k>0$ and
$a\in\F_p^*$. If $\gcd(k,p-1)=5$, then $p=11$. If $\gcd(k,p-1)=6$,
then $p\in\{7,13,19,31\}$. If $\gcd(k,p-1)=7$, then $p=29$. If
$\gcd(k,p-1)=8$, then $p=17$. Conversely, each of these possibilities
actually occurs for some $n,k,a$.
\end{cor}

There is no difficulty extending this to larger values of
$\gcd(k,p-1)$.


\section{Existence results}
\label{sec exist}

In this section we estimate the number of permutation binomials of
prescribed shapes.

\begin{thm}
\label{cw}
Let $n,k>0$ be integers with $\gcd(n,k,q-1)=1$, and suppose
$q\ge 4$.  If $\gcd(k,q-1) > 2q(\log\log q)/\log q$, then
there exists $a\in\F_q^*$ such that $x^n(x^k+a)$
permutes\/ $\F_q$.  Further, letting $T$ denote the number of $a\in\F_q$
for which $x^n(x^k+a)$ permutes\/ $\F_q$, and writing $r:=(q-1)/\gcd(k,q-1)$,
we have
\begin{align*}
\frac{r!}{r^r}&\left(q+1-\sqrt{q}(r^{r+1}-2r^r-r^{r-1}+2) - (r+1)r^{r-1}\right)
\le T \\
&\le \frac{r!}{r^r}\left(q+1+\sqrt{q}(r^{r+1}-2r^r-r^{r-1}+2)\right).
\end{align*}
\end{thm}

\begin{cor}
For fixed $r$, as $q\to\infty$ we have $T\sim q(r!)/r^r$.
\end{cor}

Note that Stirling's approximation says that $r!/r^r$ is asymptotic to
$\sqrt{2\pi r}/e^r$ as $r\to\infty$.

We will prove Theorem~\ref{cw} as a consequence of several lemmas, which we
suspect will be useful in future work improving the bounds in Theorem~\ref{cw}.
In these lemmas, $\mu_r$ denotes the set of $r^{\operatorname{th}}$ roots of
unity in $\F_q$, and $\Sym(\mu_r)$ denotes the set of permutations of $\mu_r$.

\begin{lemma}
\label{l1} Let $k,n>0$ be integers with $k\mid (q-1)$ and
$\gcd(n,k)=1$, and put $r:=(q-1)/k$.  For $a\in\F_q$, the polynomial
$f(x):=x^n(x^k+a)$ permutes $\F_q$ if and only if there exists
$\pi\in\Sym(\mu_r)$ such that every $\zeta\in\mu_r$ satisfies
$(\zeta+a)^k = \pi(\zeta)/\zeta^n$.
\end{lemma}

\begin{proof}
For $\delta\in\mu_k$ we have $f(\delta x)=\delta^n f(x)$; since $\gcd(n,k)=1$,
it follows that the values of $f$ on $\F_q$ comprise all the
$k^{\operatorname{th}}$ roots of the values of
$f(x)^k=x^{kn}(x^k+a)^k$.  Thus, $f$ permutes $\F_q$ if and only if
$g(x):=x^n (x+a)^k$ permutes the set of $k^{\operatorname{th}}$
powers in $\F_q$, or in other words $g$ permutes $\mu_r$.  Writing
$\pi$ for the map $\mu_r\to\F_q$ induced by $g$, the result
follows.
\end{proof}

Next we restate Lemma~\ref{l1} in terms of solutions to a system of
nonlinear equations over $\F_q$.  In this statement,
$\nu:\mu_r\to\F_q^*$ is a fixed map with the property that
$\nu(\zeta)^k=\zeta$ for every $\zeta\in\mu_r$.

\begin{lemma}
\label{l5} Let $k,n,r$ be as in Lemma~\ref{l1}. For $a\in\F_q$, the
polynomial $f(x):=x^n(x^k+a)$ permutes $\F_q$ if and only if there
exists $\pi\in\Sym(\mu_r)$ such that, for each $\zeta\in\mu_r$, there
is a solution $y_{\zeta}\in\F_q^*$ to the equation $\zeta+a =
y_{\zeta}^r\nu(\pi(\zeta)/\zeta^n)$. Moreover, for any fixed $a\in\F_q$,
there is at most one such permutation $\pi$.
\end{lemma}

\begin{proof}
By Lemma~\ref{l1}, $f$ permutes $\F_q$ if and only if there exists
$\pi\in\Sym(\mu_r)$ such that $(\zeta+a)^k = \pi(\zeta)/\zeta^n$ for
all $\zeta\in\mu_r$.  This equation shows that at most one $\pi$
corresponds to a given $f$.  For fixed $\pi$ and $\zeta$, the
equation is equivalent to the existence of $y_{\zeta}\in\F_q^*$ such
that $\zeta+a = y_{\zeta}^r\nu(\pi(\zeta)/\zeta^n)$.
\end{proof}

Let $A$ be transcendental over $\F_q$, and for $\pi\in\Sym(\mu_r)$
let $F_{\pi}= \F_q(\{Y_{\zeta}:\zeta\in\mu_r\})$ where
$Y_{\zeta}^r\nu(\pi(\zeta)/\zeta^n)= \zeta+A$.  We will translate
Lemma~\ref{l5} into a statement about $F_{\pi}$, which will enable
us to apply Weil's bound on the number of degree-one places of a
function field over a finite field.  In order to make this
translation, we need to know some basic facts about $F_{\pi}$, which
we record in the next lemma.  In the remainder of this section we
use various standard facts about algebraic function fields, for
which a convenient reference is~\cite{St}.

\begin{lemma}
\label{l2} Let $k,n,r$ be as in Lemma~\ref{l1}. Then $\F_q$ is
algebraically closed in $F_{\pi}$, and $F_{\pi}/ \F_q(A)$ is Galois
with group $(\Z/r\Z)^r$.  Moreover, the extension $F_{\pi}/\F_q(A)$
has ramification index $r$ over $A=\infty$ and $A\in -\mu_r$, and is
unramified over all other places of\/ $\F_q(A)$. The genus of
$F_{\pi}$ is $(r^{r+1}-2r^r-r^{r-1}+2)/2$.
\end{lemma}

\begin{proof}
Let $E_{\zeta}$ be the field $\F_q(Y_{\zeta})$.  Then
$E_{\zeta}/\F_q(A)$ is a degree-$r$ Kummer extension which is
totally ramified over $A=\infty$ and $A=-\zeta$, and unramified over
all other places.  Since each extension $E_{\zeta}/\F_q(A)$ is
totally ramified over a place which does not ramify in any other
$E_{\zeta'}/\F_q(A)$, it follows that the compositum $F_{\pi}$ of
the various fields $E_{\zeta}$ is a degree-$r^r$ extension of $\F_q(A)$
such that $\F_q$ is algebraically closed in $F_{\pi}$.  Moreover, $F_{\pi}$
is a Galois extension of $\F_q(A)$ with Galois group $(\Z/r\Z)^r$.
By Abhyankar's lemma, $F_{\pi}/\F_q(A)$ has
ramification index $r$ over $A=\infty$ and $A\in -\mu_r$, and this
extension is unramified over all other places of $\F_q(A)$.  Now the
Riemann-Hurwitz formula yields the genus of $F_{\pi}$.
\end{proof}

Now we can restate Lemma~\ref{l5} in terms of places of $F_{\pi}$:

\begin{lemma}
\label{l3}
Let $k,n,r$ be as in Lemma~\ref{l1}.
For $a\in\F_q$, the polynomial $f(x):=x^n(x^k+a)$ permutes $\F_q$
if and only if there exists $\pi\in\Sym(\mu_r)$ such that $F_{\pi}$ has
a degree-one place with $A=a$ and every $Y_{\zeta}\ne 0$.
Moreover, for any fixed $a\in\F_q$, there is at most one such
permutation $\pi$.
\end{lemma}

\begin{proof}[Proof of Theorem~\ref{cw}]
Fix $k,n,r$.  As in the proof of Theorem~\ref{gen}, we may assume $k\mid (q-1)$.
Pick a permutation $\pi\in\Sym(\mu_r)$ and a map
$\nu:\mu_r\to\F_q^*$ such that $\nu(\zeta)^k=\zeta$ for every $\zeta\in\mu_r$.
Let $N_{\pi}$ denote the number of degree-one places of $F_{\pi}$.
Then Weil's bound gives
\[
|N_{\pi} - (q+1)| \le (r^{r+1}-2r^r-r^{r-1}+2)\sqrt{q}.
\]
The ramified places in $F_{\pi}/\F_q(A)$ are precisely the places of
$F_{\pi}$ for which either $A=\infty$ or some $Y_{\zeta}\in\{0,\infty\}$.
The number of such places is at most $(r+1)r^{r-1}$.  All other rational
places of $F_{\pi}$ occur in $\Gal(\F_{\pi}/\F_q(A))$-orbits of size $r^r$,
with each orbit corresponding to a unique place of $\F_q(A)$.  Let $T$ denote
the number of values $a\in\F_q$ for which $x^n(x^k+a)$ permutes $\F_q$.
By Lemma~\ref{l3} we have
\begin{align*}
r!&\frac{q+1-(r^{r+1}-2r^r-r^{r-1}+2)\sqrt{q} - (r+1)r^{r-1}}{r^r} \le T \\
&\le r!\frac{q+1+(r^{r+1}-2r^r-r^{r-1}+2)\sqrt{q}}{r^r}.
\end{align*}
In particular, $T>1$ whenever $q>r^{2r+2}$ and $q>2$.
The former inequality is true whenever $q\ge 7$ and
$r<(\log q)/(2\log\log q)$, or equivalently $q\ge 7$ and
\[
k > \frac{2(q-1)\log\log q}{\log q}.
\]
For $q\in\{4,5\}$ we have $2q(\log\log q)/\log q
>(q-1)/2$, so it remains to show that there are permutation binomials
$x^n(x^{q-1}+a)$ (with $a\ne 0$) for every $n$ coprime to $q-1$.
By Lemma~\ref{l1}, this binomial permutes $\F_q$
whenever $a\in\F_q^*\setminus\{-1\}$.
\end{proof}

\begin{remark}
In this proof, we treated the various $\pi$'s independently.  This is
inefficient, especially since distinct $\pi$'s give disjoint sets of $a$'s.
If one could combine the information from distinct $\pi$'s more effectively,
it might be possible to remove the $\log\log q$ factor from Theorem~\ref{cw}.
We now take a first step in this direction (based on an idea in~\cite{CW}),
by effectively combining
the information from $r$ distinct $\pi$'s.  To start with, consider any of
the $(r-1)!$ permutations $\pi_0\in\Sym(\mu_r)$ with $\pi_0(1)=1$.  Now the
`$\zeta=1$' equation $(1+a)^k=\pi(1)$ can be used
as the definition of $\pi(1)$ (so long as $a\ne -1$), and we seek solutions for
each of the $(r-1)!$ permutations $\pi = (1+a)^k\cdot \pi_0$.  Thus, for each
such $\pi$, we pick $\nu$ as before and consider the function field defined by
$Y_{\zeta}^r \nu(\pi_0(\zeta)/\zeta^n)=(\zeta+A)/(1+A)$.  By the same method
as above, we find
that
\[
\frac{q-2\sqrt{q}+1}{r^{r-1}} - (r-3)\sqrt{q} - 2 \le \frac{T}{(r-1)!} \le
\frac{q+2\sqrt{q}+1}{r^{r-1}} + (r-3)\sqrt{q}.
\]
Here, as usual, one can obtain better bounds by applying the various
improvements to the Weil bound due to Manin~\cite{Manin},
Ihara~\cite{Ihara}, Drinfel'd-Vl\u{a}du\c{t}~\cite{DV},
Serre~\cite{Se1,Se2}, Oesterl\'e~\cite{Se2},
St\"ohr-Voloch~\cite{SV}, etc.

The following variant was noted implicitly in~\cite{CW} and
explicitly in~\cite{WL}: if $q$ is sufficiently
large compared to $r$ and $q\equiv 1\pmod{r}$, then there exists $a\in\F_q^*$ such
that, for every $n,k>0$ with $\gcd(n,q-1)=1$ and
$\gcd(k,q-1)=(q-1)/r$, the polynomial $x^n(x^k+a)$ permutes $\F_q$.
The novel feature here is that a single $a$ works for every $n$ and
$k$; one unfortunate aspect is that we need $\gcd(n,q-1)=1$,
whereas in Theorem~\ref{cw} we required only that
$\gcd(n,(q-1)/r)=1$. The modified proof described in this remark
gives a quantitative version of this result, so long as we restrict to
$\pi_0$ being the identity.  Let $\hat T$ denote the number of
values $a\in\F_q$ such that, for every $n,k>0$ with $\gcd(n,q-1)=1$
and $\gcd(k,q-1)=(q-1)/r$, the polynomial $x^n(x^k+a)$ permutes
$\F_q$.  Our proof in this remark (with $\pi_0(x)=x$) shows that
\[
\hat T \ge (q-2\sqrt{q}+1)/r^{r-1} - \sqrt{q}(r-3) - 2.
\]
\end{remark}

\begin{remark}
In case $r=2$, the function field $F_{\pi}$ occurring in the proof of
Theorem~\ref{cw} has genus zero, and hence can be parametrized.
This leads to explicit expressions for
the allowable values of `$a$' in this case~\cite{C,NR,W2}.
For larger values $r$, the field $F_{\pi}$ has larger genus, so one
does not expect a simple exact formula for its number of rational
places.  And indeed, already for $r=3$ the data suggests there is no
simple formula for the number of $a\in\F_q$ such that
$x(x^{(q-1)/r}+a)$ permutes $\F_q$, or more generally for the number
of permutation binomials of degree less than $q$ for which $(q-1)/r$
is the gcd of $q-1$ with the difference between the degrees of the
terms. A priori it is conceivable that there might be a nice formula
for the latter number but no nice formula for the former, since the
latter corresponds to the sum of the numbers of rational places on
the various fields $F_{\pi}$; however, the data suggests there are
no nice formulas when $r>2$.
\end{remark}

\begin{remark}
Theorem~\ref{cw} is a refinement of a result of Carlitz and
Wells~\cite{CW}. Our version differs from the original one in
various ways: it is effective, it gives an estimate on the number of
permutation binomials of prescribed shapes, it applies when
$\gcd(n,k,q-1)=1$ rather than $\gcd(n,q-1)=1$, and the proof is
geometric (in contrast to the intricate manipulation of character
sums in~\cite{CW}).  Still, we emphasize that the key idea of using
the Weil bound to prove existence of permutation binomials is due to
Carlitz~\cite{C}.
\end{remark}

\section{Heuristic}
\label{heuristic}

In this section we give a heuristic suggesting that `at random' there would
not be any permutation binomials $x^m+ax^n$ over $\F_q$ (with $m>n>0$)
such that $\gcd(m-n,q-1)<q/(2\log q)$, at least for $q$ sufficiently large.

As in the proof of Theorem~\ref{gen}, it suffices to consider
$f(x):=x^n(x^k+a)$ where $k\mid (q-1)$ and $n$ is coprime to $k$.
By Lemma~\ref{l1}, for fixed $k$, we
need only consider a single such value $n$ in each class
modulo $(q-1)/k$ which contains integers coprime to $k$.
Further, since composing $f(x)$ on both sides with scalar multiples
does not affect whether $f(x)$ permutes $\F_q$, 
we need only consider $a$'s representing the distinct cosets of the
$k^{\operatorname{th}}$ powers in $\F_q^*$ (for fixed $k$ and $n$).
Thus, for fixed $k$, there are fewer than $q$ polynomials to consider.
Since $\gcd(n,k)=1$, the values of $f$ comprise all the
$k^{\operatorname{th}}$ roots of the values of $f^k$; but the latter
are just $0$ and the values of $x^n(x+a)^k$ on $(\F_q^*)^k$. Thus,
$f$ permutes $\F_q$ if and only if $g(x):=x^n(x+a)^k$ permutes
$(\F_q^*)^k$.  Note that $(\F_q^*)^k$ equals the group $\mu_r$ of
$r^{\operatorname{th}}$ roots of unity in $\F_q^*$, where
$r:=(q-1)/k$. Here $g$ maps $\mu_r$ into $\mu_r$ if and only if
$(-a)^r\ne 1$, which we assume in what follows.  Now, the
probability that a random mapping $\mu_r\to\mu_r$ is bijective
is $r!/r^r$.  Assuming that $g$ behaves like a random map, the
expected number of permutation binomials of the form $x^n(x^k+a)$
(for fixed $q$, after our various equivalences on $n$, $k$, $a$) is
at most $q(r!)/r^r$.  Restricting to $k<q/(2\log q)$ and summing
over all $q$, we get an expected number
\[
E:=\sum_q \sum_{\substack{{r\mid (q-1)}\\ {r>2\log q}}}
q\frac{r!}{r^r}.
\]
We now show that $E$ is finite.  By reversing the order of
summation, we find that $E=\sum_{r=1}^{\infty} (r!/r^r) F(r)$, where
\[
F(r):=\sum_{\substack{{q<e^{r/2}}\\
{q\equiv 1\, ({\text{mod }} r)}\\{q \,\text{ prime power}}}} q.
\]
The number of prime powers less than $x$ which are not prime is
at most
\[
\sum_{n=2}^{\lfloor\log_2 x\rfloor} x^{1/n} < \sqrt{x}+\sqrt[3]{x}\log_2 x.\]
Thus, for fixed $r$, the number of nonprime $q$ which contribute to $F(r)$ is at most
$e^{r/4}+e^{r/6}r/(2\log 2)$.  By the Brun--Titchmarsh theorem
\cite[Thm.~3.8]{HR}, the number of prime $q$ which contribute to $F(r)$ is
at most
\[
\frac{3e^{r/2}}{\phi(r)\log \frac{e^{r/2}}{r}}.\]
Since
\[
\phi(r)> \frac{r}{e^{\gamma}\log\log r + \frac{3}{\log\log r}}
\]
for $r\ge 3$ (\cite[Thm.~15]{RS}), for $r\ge 3$ we have
\[
\frac{F(r)}{e^r}\le \frac{3(e^{\gamma}\log\log r + \frac{3}{\log\log r})}
{r(\frac{r}2-\log r)} +
\frac{1}{e^{r/4}}+\frac{r}{2e^{r/3}\log 2}.
\]
Using Stirling's inequality $r!<(r/e)^r\sqrt{2\pi r}e^{1/{12r}}$, we get
\[
E\le\sum_{r=3}^{\infty} \sqrt{2\pi r} e^{\frac{1}{12r}}
\left( \frac{3e^{\gamma}\log\log r + \frac{9}{\log\log r}}{r(\frac{r}2-\log r)} + \frac{1}{e^{r/4}} + \frac{r}{2e^{r/3}\log 2}\right),
\]
which is finite.  By combining the above bounds on $F(r)$ with explicit
calculation of the first few values of $F(r)$, we find that $E<40$.

Since $E$ is finite (and small), we expect that `at random' there
would be few (or no)
permutation binomials $x^m+ax^n$ over $\F_q$ with $m>n>0$ and
$\gcd(m-n,q-1)<q/(2\log q)$.

We used a computer to verify that, for $p<10^5$, there are no permutation
binomials $x^m+ax^n$ over $\F_p$ with $m>n>0$ and
$\gcd(m-n,p-1)<p/(2\log p)$.  Combined with the above heuristic, this leads
us to conjecture that the same conclusion holds for all primes $p$.

On the other hand, the heuristic applies to nonprime fields as well,
and for those fields we know some infinite families of
counterexamples. For instance, in~\cite{TZ}, Tom Tucker and the
second author showed that $x^{p+2}+ax$ permutes $\F_{p^2}$ whenever
$\#\langle a^{p-1}\rangle=6$. Several additional examples can be
found in~\cite{TZ}, and we will present further examples in a
forthcoming paper. However, every known counterexample over a
nonprime field $\F_q$ has unusual properties related to the
subfields of $\F_q$; thus, we view these examples as violating the
randomness hypotheses of our heuristic, rather than the heuristic
itself.


\appendix
\section*{Appendix}
In this appendix we prove the following result:

\theoremstyle{plain}
\newtheorem*{thma}{Theorem~\ref{WT}}
\begin{thma}
If $x^m+ax^n$ permutes the prime field\/ $\F_p$, where $m>n>0$ and
$a\in\F_p^*$, then $p-1 \le (m-1)\cdot\max(n,\gcd(m-n,p-1))$.
\end{thma}

As noted in the introduction, this result follows from
Theorem~\ref{intro1} in all cases except when $n=1$ and $(m-1)\mid
(p-1)$.  However, the proof we present here is quite different from
the proof of Theorem~\ref{intro1}, so the method might well be
useful in other investigations. Theorem~\ref{WT} may be viewed as
the `least common generalization' of a result of Wan and a result of
Turnwald. Our proof uses ideas from both of their proofs. Wan's
result~\cite[Thm.\ 1.3]{W} is

\newtheorem*{thmu}{Theorem}
\begin{thmu} If $x^m+ax$ permutes the prime field\/ $\F_p$, where $m>1$ and $a\in\F_p^*$,
then $p-1\le (m-1)\cdot \gcd(m-1,p-1)$.
\end{thmu}

Turnwald's result~\cite[Thm.\ 2]{T} is

\begin{thmu}
If $x^m+ax^n$ permutes $\F_p$, where
$m>n>0$ and $a\in\F_p^*$, then $p<m\cdot\max(n,m-n)$.
\end{thmu}

\begin{proof}[Proof of Theorem~\ref{WT}]
Suppose  $f(x):=x^m + ax^n$  permutes $\F_p$, where $m>n>0$ and $a\in\F_p^*$.
If $f(x) = \hat f(x^e)$, then the desired inequality for $f$
would follow from the corresponding inequality for $\hat f$; thus, we
may assume $\gcd(m,n)=1$.  Moreover, since $f$ permutes $\F_p$ we
have $\gcd(m-n,p-1)>1$ (since otherwise $f$ has more than one root),
so $n\le m-2$ and $m\ge 3$.  Write $p=mk+r$ with $0\le r<m$.  Since
$\gcd(n,m-n)=1$, there are integers $u,v$ with $nu-(m-n)v=r-1$; we may
assume $0<u\le m-n$.  Thus
\[v=(nu-r+1)/(m-n) \le n + 1/(m-n) < n+1,\] so $v\le n$.
Also $v > (n-m+1)/(m-n) > -1$, so $v\ge 0$.

If $v>k$, then (since $k=\lfloor{p/m}\rfloor$) we have
$p < mv \le mn$, so the result holds.  Henceforth we assume $v\le k$.
Moreover, since $\gcd(m-n,p-1)\ge 2$, the result is clear when $m>p/2$;
thus, we assume $m\le p/2$.  Since $3\le m$, this implies $p\ge 7$ and $m<p-3$.

We will use Hermite's criterion with exponent $k+u$.
Before doing so, we show that $0<k+u<p-1$.
The first inequality is clear, since $u>0$ and $k=\lfloor p/m\rfloor \ge 0$.
Now,
\[
k+u = \left\lfloor\frac pm\right\rfloor+u\le \frac pm+u \le\frac
pm+m-n \le \frac pm+m-1.\] Since $p>m+3$ (and $m\ge 3$), we have
$p>m^2/(m-1)$, so $m<p(m-1)/m$ and thus $p/m+m < p$.  Hence
$k+u<p-1$.

Since $0<k+u<p-1$, we have $p\nmid \binom{k+u}{t}$ for $0\le t\le
k+u$; hence the degrees of the terms of $f^{k+u}$ are precisely the
numbers $mt+n(k+u-t)$ with $0\le t\le k+u$. Since \[p-1 = mk+(r-1) =
mk + nu - (m-n)v = m(k-v) + n(u+v),\] there is a term of degree
$p-1$. Since $f$ is a permutation polynomial, Hermite's criterion
implies there must be another term of degree divisible by $p-1$.
Thus, there exists $A\ne k-v$ with $0\le A\le k+u$ such that
$mA+n(k+u-A)\equiv 0\bmod{(p-1)}$. Since increasing $t$ will increase
the value of $mt+n(k+u-t)$, and the value of this quantity for $t=A$
is larger than the corresponding value for $t=k-v$, it follows that
$A>k-v$.  Subtracting, we get $m(A-(k-v)) + n(k-A-v)\equiv
0\bmod{(p-1)}$, so $p-1$ divides $(m-n)(A-(k-v))$. In other words,
$(p-1)/\gcd(p-1,m-n)$ divides $A-(k-v)$. Since $A>k-v$, this implies
\[
\frac{p-1}{\gcd(p-1,m-n)} \le A-(k-v) \le (k+u)-(k-v) = u+v.\] Since
$u\le m-n$ and $v\le n$, we have $u+v\le m$; however, equality
cannot hold, since it would imply that $r-1=nu-(m-n)v=0$ so $r=1$,
whence $p-1=p-r=mk$, which is a contradiction since $m>1$ is the
degree of a permutation polynomial.  Thus $u+v\le m-1$, so $p-1 \le
(m-1)\cdot\gcd(p-1,m-n)$.
\end{proof}

\end{document}